\documentclass[titlepage,12pt]{article}

\usepackage{amssymb}
\usepackage{latexsym}

\newcommand{\ben}{\begin{enumerate}}
\newcommand{\een}{\end{enumerate}}
\newcommand{\ble}{\begin{lem}}
\newcommand{\ele}{\end{lem}}
\newcommand{\bth}{\begin{thm}}
\renewcommand{\eth}{\end{thm}}
\newcommand{\bpr}{\begin{prop}}
\newcommand{\epr}{\end{prop}}
\newcommand{\bco}{\begin{cor}}
\newcommand{\eco}{\end{cor}}
\newcommand{\bcon}{\begin{conj}}
\newcommand{\econ}{\end{conj}}
\newcommand{\bde}{\begin{defn}}
\newcommand{\ede}{\end{defn}}
\newcommand{\bex}{\begin{exa}}
\newcommand{\eex}{\end{exa}}
\newcommand{\barr}{\begin{array}}
\newcommand{\earr}{\end{array}}
\newcommand{\btab}{\begin{tabular}}
\newcommand{\etab}{\end{tabular}}
\newcommand{\beq}{\begin{equation}}
\newcommand{\eeq}{\end{equation}}
\newcommand{\bea}{\begin{eqnarray*}}
\newcommand{\eea}{\end{eqnarray*}}
\newcommand{\bce}{\begin{center}}
\newcommand{\ece}{\end{center}}
\newcommand{\bpi}{\begin{picture}}
\newcommand{\epi}{\end{picture}}
\newcommand{\bfi}{\begin{figure} \begin{center}}
\newcommand{\efi}{\end{center} \end{figure}}
\newcommand{\capt}{\caption}
\newcommand{\bsl}{\begin{slide}{}}
\newcommand{\esl}{\end{slide}}

\newcommand{\bib}{thebibliography}
\newcommand{\pf}{{\bf Proof.}}
\newcommand{\qed}{\rule{1ex}{1ex}}
\newcommand{\Qed}{\rule{1ex}{1ex} \medskip}

\newcommand{\sbs}{\subset}
\newcommand{\sbe}{\subseteq}

\newcommand{\setm}{\setminus}

\newcommand{\mh}{\hat{0}}
\newcommand{\Mh}{\hat{1}}

\newcommand{\ce}[1]{\lceil #1 \rceil}

\newcommand{\fall}[2]{\langle{#1}\rangle_{#2}}
\newcommand{\ffall}[2]{\langle\langle{#1}\rangle\rangle_{#2}}

\newcommand{\al}{\alpha}

\newcommand{\de}{\delta}
\newcommand{\ep}{\epsilon}

\newcommand{\bZ}{{\bf Z}}
\newcommand{\bbC}{{\mathbb C}}
\newcommand{\bbF}{{\mathbb F}}

\newcommand{\bbR}{{\mathbb R}}

\newcommand{\cA}{{\cal A}}
\newcommand{\cB}{{\cal B}}

\newcommand{\cD}{{\cal D}}

\newcommand{\cH}{{\cal H}}

\newcommand{\cS}{{\cal S}}

\newcommand{\cW}{{\cal W}}

\newcommand{\alt}{\tilde{\alpha}}

\newcommand{\aim}{Adv. in Math.}

\newcommand{\dm}{Discrete Math.}

\newcommand{\im}{Invent. Math.} 

\newcommand{\jac}{J. Algebraic Combin.}

\newcommand{\jams}{J. Amer. Math. Soc.}

\newtheorem{thm}{Theorem}[section]
\newtheorem{prop}[thm]{Proposition}
\newtheorem{cor}[thm]{Corollary}
\newtheorem{lem}[thm]{Lemma}
\newtheorem{conj}[thm]{Conjecture}
\newtheorem{exa}[thm]{Example}
\newtheorem{defn}[thm]{Definition}

\begin{document}
\pagestyle{empty}
\title{Characteristic and Ehrhart Polynomials}
\author{Andreas Blass\\
Department of Mathematics \\ University of Michigan \\
Ann Arbor, MI 48109-1003\\
ablass@umich.edu\\[5pt]
and\\[5pt]
Bruce E. Sagan  \\
Department of Mathematics \\ Michigan State University
\\ East Lansing, MI 48824-1027\\sagan@math.msu.edu}

\date{\today \\[1in]
        \begin{flushleft}
        \end{flushleft}
       }
\maketitle

\begin{flushleft} Proposed running head: \end{flushleft}
        \begin{center} 
Characteristic and Ehrhart Polynomials
        \end{center}

Send proofs to:
\begin{center}
Bruce E. Sagan \\ Department of Mathematics \\Michigan State
University \\ East Lansing, MI 48824-1027\\[5pt]
Tel.: 517-355-8329\\
FAX: 517-336-1562\\
Email: sagan@math.msu.edu
\end{center}

        \begin{abstract} 
Let $\cA$ be  a subspace arrangement and let
$\chi(\cA,t)$ be the characteristic polynomial of its intersection
lattice $L(\cA)$.
We show  that if the subspaces in $\cA$
are taken from $L(\cB_n)$, where $\cB_n$ is the type $B$ Weyl
arrangement, then $\chi(\cA,t)$ counts a certain set of 
lattice points.  One can use  this result to study the partial
factorization of $\chi(\cA,t)$ over the integers and the coefficients
of its expansion in various bases for the polynomial ring $\bbR[t]$.
Next we prove that the characteristic
polynomial of any Weyl hyperplane arrangement can be expressed in
terms of an Ehrhart quasi-polynomial for its affine Weyl chamber.
Note that our first result deals with all subspace
arrangements embedded in $\cB_n$ while the second deals  with
all finite Weyl groups but only their hyperplane arrangements.
\end{abstract}
\pagestyle{plain}

\section{Introduction and Background}

An {\it arrangement\/} is a finite set 
        \beq                                            \label{AH1}
        \cA=\{K_1, \ldots, K_m\}
        \eeq
of proper subspaces of Euclidean space $\bbR^n$.  All the subspaces we
consider will be linear and so go through the origin.
If each $K_i$ has
dimension $n-1$, then $\cA$ is called a {\it hyperplane
arrangement\/}.  We sometimes refer to general arrangements as {\it
subspace arrangements\/} to emphasize that they need not be hyperplane
arrangements.  We write $\bigcup\cA$ for the set-theoretic union of
the subspaces in $\cA$, i.e., $\bigcup_{i=1}^m K_i$.

The theory of hyperplane arrangements is a beautiful area of
mathematics which brings together ideas from topology, algebra, and
combinatorics.  Its roots go back to the end of the 19th century but
it is also an active area of research today.  The recent book
\cite{ot:ah} of Orlik and Terao covers both classical work and recent
developments in the field.  Subspace arrangements, on the other hand,
have received relatively little attention yet, as was noted in the recent
survey article of Bj\"orner~\cite{bjo:sa}.  It is important to
emphasize that in most cases it is {\it not} easy to generalize
results from the hyperplane case to the subspace case.  Particularly
nicely behaved hyperplane arrangements are those which are associated
with finite Weyl groups (see, e.g., \cite{os:ca}).  We wish to
study these arrangements and certain subspace arrangements related to
them.  We begin by establishing some notation and terminology.

Let $\cA$ be an arrangement as in (\ref{AH1}) above, and assume, for
simplicity, that there are no containments among the $K_i$.  Let
$L=L(\cA)$ be the set of all intersections of these subspaces, ordered
by reverse inclusion, called the {\it intersection lattice}.
(Concepts from lattice theory that are not explained here can be found
in Stanley's text~\cite{sta:ec1}.)  Note that $L$ has a unique minimal
element $\mh$ corresponding to $\bbR^n$, an atom corresponding to each
$K_i$, and a unique maximal element $\Mh$ corresponding to
$\bigcap_{i=1}^m K_i$.  If $\cA$ is a hyperplane arrangement then
$L(\cA)$ is a geometric lattice, but in general it is not even ranked.
If $\cA$ and $\cB$ are subspace arrangements such that $\cA\sbe
L(\cB)$, i.e. all the subspaces in $\cA$ are intersections of
subspaces in $\cB$, then we say that $\cA$ is {\it embedded} in $\cB$.

Given an arrangement $\cA$,
let $\mu(X)=\mu(\mh,X)$ denote the {\it M\"obius
function} of the lattice $L(\cA)$; it is uniquely defined by
        $$\sum_{Y\le X} \mu(Y)=\de_{\mh,X}$$
where $\de_{\mh,X}$ is the Kronecker delta.
The M\"obius function is one of the fundamental invariants of any
partially ordered set; see the seminal article of Rota~\cite{rot:tmf}.
The {\it characteristic polynomial\/} of $\cA$ is
        \beq                                                    \label{chi}
        \chi(\cA,t)=\sum_{X\in L(\cA)}\mu(X)t^{\dim X}.
        \eeq
Since the characteristic
polynomial is just the generating function for the M\"obius function,
it is also of prime importance.
Our results in this paper give a combinatorial interpretation for the
characteristic polynomials of hyperplane arrangements associated to
Weyl groups and subspace arrangements embedded in some of these
Weyl arrangements. 

For any finite Weyl group, $W$, there is a corresponding hyperplane
arrangement $\cW$ whose elements are the reflecting hyperplanes of
$W$.  Initially we shall be interested in the case where $W$ comes from
one of the three infinite families $A_n, B_n, D_n$.  (The arrangement
for $C_n$ is clearly the same as that for $B_n$.)  In terms of the
coordinate functions $x_1,\ldots, x_n$ in $\bbR^n$, the associated
hyperplane arrangements can be defined as 
        \bea 
        \cA_n&=&\{x_i=x_j\ :\ 1\le i<j\le n\},\\ 
        \cD_n&=&\cA_n\cup\{x_i=-x_j\ :\ 1\le i<j\le n\},\\
        \cB_n&=&\cD_n\cup\{x_i=0\ :\ 1\le i\le n\} 
        \eea 
so that $\cA_n\sbs\cD_n\sbs\cB_n$.  Note that $n$ here refers to the
dimension of the space, not the number of fundamental reflections
(which is $n-1$ for $\cA_n$ and $n$ for the other two).

\section{Arrangements Embedded in $\cB_n$}

We shall now give our first main result: a combinatorial
interpretation for the characteristic 
polynomial of any subspace arrangement embedded in one of the three
infinite families of Weyl hyperplane arrangements.  It was obtained
in an attempt to generalize Zaslavsky's beautiful theory of signed
graph coloring~\cite{zas:grs,zas:sgc,zas:cis}.  Given integers $r\le
s$, we let $[r,s]=\{r, r+1,\ldots, s\}$.  Note that if $r=-s$ then
$t=|[-s,s]|$ is odd, where $|\cdot|$ denotes cardinality.
Note also that $[-s,s]^n$ is just the cube of
points in $\bZ^n$ centered at the origin with $t$ points on a side.
So $[-s,s]^n\setm\bigcup\cA$ is the set of points of $\bZ^n$
that are in this cube but not on any subspace from $\cA$.

\bth            \label{thb}
If $\cA\sbe L(\cB_n)$ then for any $t=2s+1$
        $$\chi(\cA,t)=|[-s,s]^n\setm\bigcup\cA|.$$
\eth

Note that the hypothesis of the theorem does not preclude the
possibility that $\cA$ may also be embedded in $\cA_n$ or $\cD_n$, as
these are embedded in $\cB_n$.  Let us give a concrete example of
this result before proving it.  Let 
        $$\cA=\cB_2=\{x=0,y=0,x=y,x=-y\}.$$ 
Also let  $s=2$ so that $t=5$.  
Then $[-2,2]^2$ and $\cB_2$ are shown in Figure~\ref{Bfig}.  Removing the
lines of $\cB_2$ from the cube leaves 8 lattice points.  On the
other hand it is well known that $\chi(\cB_2,t)=(t-1)(t-3)$; see equation~(\ref{Beq}).
So $\chi(\cB_2,5)=4\cdot 2=8$ as expected.

\thicklines
\setlength{\unitlength}{2pt}
\bfi
\bpi(60,30)(-10,0)
\put(0,0){\circle*{3}}
\put(10,0){\circle*{3}}
\put(20,0){\circle*{3}}
\put(30,0){\circle*{3}}
\put(40,0){\circle*{3}}
\put(0,10){\circle*{3}}
\put(10,10){\circle*{3}}
\put(20,10){\circle*{3}}
\put(30,10){\circle*{3}}
\put(40,10){\circle*{3}}
\put(0,20){\circle*{3}}
\put(10,20){\circle*{3}}
\put(20,20){\circle*{3}}
\put(30,20){\circle*{3}}
\put(40,20){\circle*{3}}
\put(0,30){\circle*{3}}
\put(10,30){\circle*{3}}
\put(20,30){\circle*{3}}
\put(30,30){\circle*{3}}
\put(40,30){\circle*{3}}
\put(0,40){\circle*{3}}
\put(10,40){\circle*{3}}
\put(20,40){\circle*{3}}
\put(30,40){\circle*{3}}
\put(40,40){\circle*{3}}
\put(0,0){\line(1,1){40}}
\put(20,0){\line(0,1){40}}
\put(0,20){\line(1,0){40}}
\put(0,40){\line(1,-1){40}}
\epi
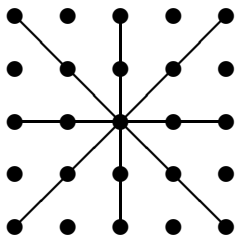
\capt{The lattice points of $[-2,2]^2\setm\bigcup\cB_2$}  \label{Bfig}
\efi

{\bf Proof of Theorem \ref{thb}.}  
 We construct two functions $f,g:L(\cA)\longrightarrow\bZ$ by
defining for each $X\in L(\cA)$
        \bea
        f(X)&=&|X\cap[-s,s]^n|,\\
        g(X)&=&|(X\setm\bigcup_{Y>X}Y)\cap[-s,s]^n|.
        \eea
Recall that $L(\cA)$ is ordered by {\it reverse} inclusion so that
$\bigcup_{Y>X}Y\sbs X$.  In particular $g(\bbR^n)=|[-s,s]^n\setm\bigcup\cA|$.   
Note also that $X\cap[-s,s]^n$ is combinatorially
just a cube of dimension $\dim X$ and side $t$ so that $f(X)=t^{\dim X}$.
Finally, $f(X) =\sum_{Y\ge X}g(Y)$ so by the M\"obius Inversion
Theorem~\cite{rot:tmf} 
        \bea
        |[-s,s]^n\setm\bigcup\cA|&=&g(\mh)\\
                &=&\sum_{X\in L(\cA)}\mu(X)f(X)\\
                &=&\sum_{X\in L(\cA)}\mu(X)t^{\dim X}\\
                &=&\chi(\cA,t)
        \eea
which is the desired result.\hfill\Qed

In the proof of Theorem~\ref{thb}, it was crucial that each of the
subspaces $X$ under consideration had exactly $t^{\dim(X)}$ points in
$[-s,s]^n$.  In fact, the {\it only} subspaces of $\bbR^n$  with this
property are those in $L(\cB_n)$.  So the method of proof of
Theorem~\ref{thb} cannot be applied directly to other arrangements.

We should also mention how our theorem is related to Zaslavsky's theory of
signed graphs.  Zaslavsky assigns to each hyperplane arrangement
$\cA$ contained (as a subset) in $\cB_n$ a signed graph $G_{\cA}$.
The graph has vertices $1,2,\ldots,n$ with a positive
(respectively, negative) edge
from vertex $i$ to vertex $j$ iff $x_i=x_j$ (respectively, $x_i=-x_j$)
is in $\cA$.  The graph $G_{\cA}$ also has a half-edge at vertex $i$ iff $x_i=0$ is in $\cA$.
He then defines a chromatic polynomial $P(G,t)$ for signed graphs 
(generalizing the one for ordinary graphs) 
and shows that $P(G_{\cA},t)=\chi(\cA,t)$.
If one thinks of the vertices of $G_{\cA}$ as being coordinates,
then a proper coloring of $G_{\cA}$ in Zaslavsky's sense turns out to
be just an element of $[-s,s]^n\setm\bigcup\cA$.  The advantages of our
viewpoint are that it applies to subspace arrangements embedded in
$\cB_n$ (not just hyperplane embeddings) and that it
admits an analog for all Weyl hyperplane arrangements
as we shall see in our second main theorem.
We should mention that Stanley~\cite{sta:gcr} has independently
formulated a version 
of Theorem~\ref{thb} for arrangements embedded in $\cA_n$ using
hypergraphs and symmetric functions.

\section{Examples}

First, let us show how Theorem~\ref{thb}  can be used to compute the
well-known
characteristic polynomials for the three infinite families of Weyl
hyperplane arrangements.  In the type $A$ case we see that a point of 
$[-s,s]^n\setm\bigcup\cA_n$ must have all coordinates different.
So there are $t=2s+1$ choices for the first coordinate, $t-1$ for the
second, etc. This gives  a total of 
        $$\chi(\cA_n,t)=t(t-1)\cdots(t-n+1).$$
It will be useful to have a notation for this falling factorial, so
we will let $\fall{t}{n}=t(t-1)\cdots(t-n+1)$.

For $\cB_n$ the points
in the cube minus the arrangement must all have different absolute
values and must be nonzero.
The first coordinate can be chosen in $t-1$ ways since zero is not
allowed.  The second coordinate can
be anything except zero and plus or minus the value of the first,
giving $t-3$ possibilities.  Continuing in this way we see that
        \beq                                    \label{Beq}
        \chi(\cB_n,t)=(t-1)(t-3)\cdots(t-2n+1).
        \eeq
We will let $\ffall{t}{n}=t(t-2)\cdots(t-2n+2)$ so that
$\chi(\cB_n,t)=\ffall{t-1}{n}$.

For the third family, note that any point of $[-s,s]^n\setm\bigcup\cD_n$
can have at most one zero coordinate.  The points with no zero
coordinate were counted in the $\cB_n$ case.  For those with one zero,
there are $n$ ways to pick this  coordinate and the remaining nonzero
ones are accounted for as in $\cB_{n-1}$.  The total is thus
$$\chi(\cD_n,t)=\chi(\cB_n,t)+n\chi(\cB_{n-1},t)=
(t-1)(t-3)\cdots(t-2n+3)(t-n+1).$$

Notice that in all three of these examples $\chi$
factors over the integers.  In fact for any Weyl
hyperplane arrangement it is well known that the roots are just the
exponents of the corresponding group~\cite{ter:gef}.  The
characteristic polynomial 
of a subspace arrangement $\cS_n$  embedded in a Weyl hyperplane
arrangement $\cH_n$ from one of the three infinite families
does not always have integral roots.
But it can happen that it factors partially and is in fact divisible
by the polynomial for a hyperplane arrangement $\cH_m$, $m\le n$.
Further, when one expands $\chi(\cS_n,t)$ in terms of the basis
$\{\chi(\cH_j,t)\ :\ j\ge0\}$ for $\bbR[t]$ the coefficients 
vanish for small $j$, thus explaining
the divisibility relation since for type $A$ and $B$ we have
$\chi(\cH_j,t)|\chi(\cH_{j+1},t)$.  Finally, the coefficients in the
basis expansion turn out to be nonnegative integers 
having a nice combinatorial interpretation which makes it
obvious when they are zero.  The next few results  will illustrate this
point.  Other examples can be found in~\cite{bs:sat,zha:sba} and are being
pursued by Sagan.

To describe the subspace arrangements that we will consider, it is
convenient to have some notation.  Let $[n]=\{1,\ldots,n\}$.  If
$I=\{i,j,\ldots,k\}\sbe[n]$ then let $x_I$ stand for the
equation
$x_i=x_j=\cdots=x_k$.  So $x_I=0$ is the system of equations $x_i=0$ for
all $i\in I$. 
Also let $\pm x_I$  represent 
the set of all equations  of the form
        $$\ep_i x_i=\cdots=\ep_k x_k$$
for $\ep_i,\ldots,\ep_k\in\{\pm1\}$.  In each case we use the same
symbol to denote the corresponding subspace(s).
The {\it $k$-equal} and {\it $k,h$-equal} subspace arrangements are
defined by
        \bea
        \cA_{n,k}&=&\{x_I\ :\ I\sbe[n] \mbox{ and } |I|=k\},\\
        \cD_{n,k}&=&\{\pm x_I\ :\ I\sbe[n] \mbox{ and } |I|=k\},\\
        \cB_{n,k,h}&=&\cD_{n,k}\cup\{x_J=0\ :\ J\sbe[n] \mbox{ and } |J|=h\}.
        \eea
The $\cA_{n,k}$ arrangement first appeared in the work
of Bj\"orner, Lov\'asz and Yao~\cite{bly:ldt}, motivated by its
relevance to a certain problem in computational complexity.
Its study has been continued by these authors and Linusson,
Sundaram, Wachs and Welker in various
combinations~\cite{bl:ldt,bw:hkm,bw:snc1,bw:snc2,lin:prb,sw:hrk,sw:gal}.
The $\cB_{n,k,h}$ and $\cD_{n,k}$ were introduced by Bj\"orner and
Sagan in a paper~\cite{bs:sat} about their combinatorial and homological
properties.   Note that each of these subspace arrangements is
embedded in the hyperplane arrangement of the corresponding type and
therefore in $\cB_n$.

Consider the $k$-equal arrangement $\cA_{n,k}$ embedded in $\cA_n$
with $\chi(\cA_n)=\fall{t}{n}$.  It will be convenient
to let $S_k(n,j)$ denote the number of partitions of an $n$-element
set into $j$ subsets each of which is of size at most $k$.  Thus these
are generalizations of the Stirling numbers of the second kind.
\bth
We have the expansion
        \beq                            \label{Seq}
        \chi(\cA_{n,k},t)=\sum_j S_{k-1}(n,j) \fall{t}{j}
        \eeq
and the divisibility relation
        \beq                            \label{div}
        \fall{t}{\ce{n/(k-1)}}\ |\ \chi(\cA_{n,k},t).
        \eeq
\eth
\pf\  To get the expansion, consider an arbitrary point
$x\in[-s,s]^n\setm\bigcup\cA_{n,k}$.  So $x$ can have at
most $k-1$ of its coordinates equal.  Consider the $x$'s with exactly
$j$ different coordinates.  Then there are $S_{k-1}(n,j)$ ways to
distribute the $j$ values among the $n$ coordinates with at most $k-1$
equal.  We can then choose which values to use in
$\fall{t}{j}$ ways.  Summing over all $j$ gives the desired equation.

For the divisibility result, note that $S_{k-1}(n,j)=0$ if
$j<\ce{n/(k-1)}$ because $j$ sets of at most $k-1$ objects can
partition a set of size of at most $n=j(k-1)$.  Plugging this
into~(\ref{Seq}) finishes the proof.\hfill\Qed

We should note that expansion~(\ref{Seq}) was  derived by
Bj\"orner and Lov\'asz~\cite{bl:ldt} and by
Sundaram~\cite{sun:aht} using formal power series techniques.  Analogs of
this expansion for type $B$ and $D$ can be found in a paper of
Bj\"orner and Sagan~\cite{bs:sat} while applications to the Boolean
algebra are in Zhang's thesis~\cite{zha:sba}.

\bco Let $\cA$ be a subspace arrangement.

(a)  If $\cA$ is embedded in $\cA_n$ and we write
        \beq                            \label{ajeq}
        \chi(\cA,t)=\sum_{j=0}^n a_j \fall{t}{j}
        \eeq
then $a_j\in\bZ_{\ge0}$ for all $j$, $0\le j\le n$.  
Furthermore if $m$ is the
largest index such that $a_m=0$ then
        $$\fall{t}{m+1}\ |\ \chi(\cA,t).$$

(b)   If $\cA$ is embedded in $\cB_n$ and we write
        $$\chi(\cA,t)=\sum_{j=0}^n b_j \ffall{t-1}{j}$$
then $b_j\in\bZ_{\ge0}$ for all $j$, $1\le j\le n$.  Furthermore if
$m$ is the 
largest index such that $b_m=0$ then
        $$\ffall{t-1}{m+1}\ |\ \chi(\cA,t).$$
\eco
\pf\
We will do part (a) as (b) is similar.  Consider any $X\in L(\cA_n)$ and
define $X^0=(X\setm\bigcup_{Y>X}Y)\cap[-s,s]^n$ where $Y\in L(\cA_n)$.
Then we have 
$X^0\sbe\bigcup\cA$ if $X\sbe K$ for some $K\in\cA$.
On the other hand we have $X^0\sbe[-s,s]^n\setm\bigcup\cA$
if there is no such $K$ containing $X$.  It follows that
        $$[-s,s]^n\setm\bigcup\cA=\biguplus_{X} X^0$$
where the disjoint union is over all $X$ not contained in any subspace
of $\cA$.
Taking cardinalities on both side of this equation and using the fact
that $|X^0|=\fall{t}{\dim X}$ shows that the $a_j$ in~(\ref{ajeq}) are
nonnegative integers. 

For the divisibility relation, it suffices to prove that $a_j=0$
implies $a_{j-1}=0$.  But $a_j=0$ implies that every $X\in L(\cA_n)$
of dimension $j$ is contained in some $K\in\cA$.  Thus any $Y>X$ is in
a $K$ and $a_{j-1}=0$.\hfill\Qed

\section{Weyl Hyperplane Arrangements}

In this section we confine our attention to hyperplane
arrangements that consist of the reflecting hyperplanes of a Weyl
group.  For background information on Weyl groups, including any
concepts that we use without explanation, see the book of 
Humphreys~\cite{jh:rgcg}, whose notation we endeavor to follow.  We
shall obtain 
a combinatorial characterization of the characteristic polynomial of
such an arrangement. In rough outline, the characterization is similar
to Theorem~\ref{thb}, but the lattice $\bZ^n$ will be replaced with
another lattice, the cube of side $2s+1$ will be replaced with another
polytope, and the restriction to odd values of $t$ will be replaced
with other congruences imposed on $t$.  

Unfortunately, both of the (mathematical) meanings of ``lattice'' ---
a poset in which finite subsets have joins and meets, and a discrete
subgroup of $\bbR^n$ --- are relevant to the present discussion.  We
rely on the context to make it clear which is meant.

Let $W$ be a finite Weyl group, determined by a root system $\Phi$
spanning $\bbR^n$.  The hyperplanes orthogonal to the roots constitute
the {\it Weyl arrangement} $\cW$ associated to $W$, and the
reflections in these hyperplanes generate $W$.  Throughout this
section, we follow the convention of naming a Weyl arrangement by
the script letter corresponding to the name of the Weyl group.
This agrees with the notation in the preceding sections for $\cB_n$
and $\cD_n$, but what we now call $\cA_n$ is the restriction, to the
hyperplane $x_1+x_2+\cdots+x_{n+1}=0$, of what was previously called
$\cA_{n+1}$.  

Let $Z(\Phi)$ be the lattice in $\bbR^n$ consisting of those vectors
$x$ that satisfy $(\alpha,x)\in\bZ$ for all roots $\alpha\in\Phi$.
This is the coweight lattice associated to $\Phi$, and it will play
the role that $\bZ^n$ played in Theorem~\ref{thb}.

Our analog of the cube $[-s,s]^n$ of lattice points is
$$
P_t(\Phi)=\{x\in Z(\Phi) \mid (\alpha,x)<t \mbox{ for all $\al\in\Phi$} \}.
$$
Of course we will be interested in counting the lattice points in 
$P_t(\Phi)\setm\bigcup\cW$.

Fix a simple system 
$$
\Delta=\{\sigma_1,\dots,\sigma_n\}
$$ 
in $\Phi$.  Thus, $\Delta$ is a basis for
the vector space $\bbR^n$, and, when any root $\lambda\in\Phi$ is
written as a linear combination,
$$
\lambda=\sum_{i=1}^n c_i(\lambda)\sigma_i,
$$
of $\Delta$, the coefficients $c_i(\lambda)$ are integers and are
either all $\geq0$ or all $\leq0$.  The fact that the coefficients are
integers implies that, if a vector $x$ satisfies $(\alpha,x)\in\bZ$ for all
$\alpha\in\Delta$, then it automatically satisfies the same for all
$\alpha\in\Phi$ and therefore belongs to $Z(\Phi)$.  In other words,
in defining the coweight lattice, we could have restricted attention
to simple roots.  

If $\Phi$ is irreducible then
among all the roots there is a {\it highest} one, $\alt$,
characterized by the fact that, for all roots $\lambda$ and all
$i\in[n]$, $c_i(\alt)\geq c_i(\lambda)$.  We shall write simply $c_i$
for $c_i(\alt)$.  One final ingredient for our theorem is the {\it
index of connection}, $f$, which we define for irreducible root
systems as
\beq                    \label{fdef}
f=\frac{|W|}{n!\cdot c_1\cdots c_n}.
\eeq
For an arbitrary root system, $f$ is defined as the product of the indices
of connection for each irreducible component.
(Humphreys defines $f$~\cite[p.\ 40]{jh:rgcg} as the index of the
coroot lattice as a subgroup of the coweight lattice and derives~(\ref{fdef})
as his Proposition~4.9.  Since this formula is all we
need to know about $f$, we take it as the definition.)  

\bth                    \label{thc}
Let $\Phi$ be a root system for a finite Weyl group with associated
arrangement $\cW$.  Let $t$ be a positive integer relatively prime
to all the coefficients $c_i=c_i(\alt)$.  Then 
        $$\chi(\cW,t)=\frac{1}{f}\left|P_t(\Phi)\setm\bigcup\cW\right|.$$
\eth

\pf\
We may as well assume that $\Phi$ is irreducible since if it is not
then both sides of the given equation decompose into a product of
factors, one for each of the irreducible components.
We begin by representing vectors in a form convenient for counting the
points in $P_t(\Phi)\setminus\bigcup\cW$.  For any
$x\in\bbR^n$, let $x^*$ be the $n$-tuple consisting of the inner
products of $x$ with the simple roots, i.e., $x^*_i=(\sigma_i,x)$.  So
$x\in Z(\Phi)$ if and only if $x^*\in\bZ^n$.  Also, $x$ lies in the
open fundamental chamber $C$ of $W$ if and only if $x^*$ lies in the
open positive orthant $(\bbR_{>0})^n$.  

Since  $P_t(\Phi)$ and $\cW$ are both invariant under the
action of the group $W$, we can count the points of 
$P_t(\Phi)\setminus\bigcup\cW$ by first counting the ones in $C$ and
then multiplying by the number of chambers (which equals the group's
order $|W|$).  To do the counting in $C$, we count instead the
corresponding points $x^*$ in the positive orthant of $\bZ^n$ subject
to the requirement $x\in P_t(\Phi)$.  Note that since $x^*$ is in the
{\it open} positive orthant $x$ is automatically not in $\bigcup\cW$.
For $x^*$ in $\bZ^n$ the requirement that $x\in P_t(\Phi)$ is
equivalent to the fact that, for all roots $\lambda$, 
$$
t>(\lambda,x)=\sum_i c_i(\lambda) x^*_i.
$$
But since the $x^*_i$ are all positive, these inequalities for all
$\lambda\in\Phi$ follow from the one with the largest coefficients,
namely the one for $\lambda=\alt$.  So our task is to count the number
$\psi(t)$ of points $x^*\in(\bZ_{>0})^n$ that satisfy the one linear
inequality $\sum c_i x^*_i <t$.  This $\psi(t)$ is known as the {\it
Ehrhart quasi-polynomial} of the open simplex bounded by the
coordinate hyperplanes and the hyperplane $\sum_i c_ix^*_i=1$; see
\cite{sta:ec1}, page 235ff.
It is also interesting to note that $P_1(\Phi)$ is just the
fundamental chamber for the affine Weyl group corresponding to $W$.

Getting back to the task at hand, 
we must prove that
$\psi(t)\cdot|W|=f\cdot\chi(\cW,t)$ when $t$ is relatively prime to
all $c_i$.  Using our definition~(\ref{fdef}) of $f$ we see that this
is equivalent to showing
        $$\chi(\cW,t)=\psi(t)\cdot n!\prod_i c_i$$
for the appropriate values of $t$ and this is the form that we shall
use in practice.

To compute $\psi(t)$, we use its generating function
$\gamma(z)=\sum_t\psi(t)\cdot z^t$.  It is easy to see that the
generating function for $n$-tuples $x^*$ of positive integers with
$\sum c_i x^*_i$ equal to $t$ is 
$$
\prod_{i=1}^n(z^{c_i}+z^{2c_i}+\cdots)=
\prod_{i=1}^n\frac{z^{c_i}}{1-z^{c_i}}.
$$
To get the generating function for $\sum c_i x^*_i$ strictly smaller
than $t$, one just multiplies this by $z+z^2+z^3+\cdots$, obtaining
$$
\gamma(z)=\frac{z}{1-z}\cdot\prod_{i=1}^n\frac{z^{c_i}}{1-z^{c_i}}.
$$

If we let $m$ be the least common multiple of the $c_i$'s, then all
the fractions in this product can be written with denominator $1-z^m$.
It follows, by the general theory of rational generating functions
(cf.\ \cite{sta:ec1}, Chapter~4), that $\psi(t)$ is, for positive $t$,
a quasi-polynomial with quasi-period $m$ and degree $n$.  This means
that, when restricted to values of $t$ in any one congruence class
modulo $m$, $\psi$ is a polynomial of degree $n$.

{}From here on, the proof is computational.  One inserts into the
formula for $\gamma(z)$ the coefficients $c_i$ appropriate for a
particular $\Phi$ (cf. page~98 of \cite{jh:rgcg}), one obtains a
polynomial formula for $\psi$ on each congruence class modulo $m$
(either by direct calculation or by computing enough values of $\psi$
to uniquely interpolate polynomials of the right degree), and
one verifies that, for the congruence classes prime to $m$ (or
equivalently prime to all the $c_i$), the polynomial so obtained, when
multiplied by $|W|/f$, yields the (known) characteristic polynomial of
$\cW$.  Here are some of the computations.  

For $A_n$, the $c_i$ are all 1, so 
        $$\gamma(z)=\frac{z^{n+1}}{(1-z)^{n+1}}.$$
Here the coefficients of the expansion are well known, and we find that
$\psi(t)={t-1\choose n}$.  Multiplying by $n!\prod_i c_i=n!$ we get $\fall{t-1}n$, the
characteristic polynomial of $\cA_n$.  (This differs from the
characteristic polynomial of $\cA_n$ in the preceding section because
what was there called $\cA_n$ is the current $\cA_{n-1}$ with all
dimensions increased by 1.)

For $B_n$, the $c_i$ are all 2 except for a single 1, so $t$ is odd.
The generating function is 
$$
\gamma(z)=\frac
z{1-z}\cdot\left(\frac{z^2}{1-z^2}\right)^{n-1}\cdot\frac z{1-z} =
\frac{z^{2n}(1+z)^2}{(1-z^2)^{n+1}}.
$$
Here the expansion of $(1-z^2)^{-n-1}$
contains every even power $z^{2k}$ of $z$ with coefficient $k+n\choose
n$ (and of course contains no odd powers of $z$).  So, since $t$ is
odd, the coefficient of $z^t$ in $\gamma(z)$ is
$$
\psi(t)=2\cdot{(t-1)/2\choose n}.
$$
Multiplying by $n!\prod_i c_i=2^{n-1}n!$, we get 
$$
2^n\cdot\fall{(t-1)/2}n=\ffall{t-1}{n},
$$
the characteristic polynomial of $\cB_n$.
We note that when $t$ is even a similar
calculation gives 
$$
2^{n-1}n!\psi(t)=(t-2)(t-4)\cdots(t-2n+2)\cdot(t-n)=\chi(\cD_n,t-1).
$$
We do not know any reason for this coincidence.

The computations for $C_n, D_n$ and the exceptional root
systems follow the same pattern as those for $A_n$ and $B_n$.  The
necessary information can be found in the following table.  
In it, the $c_i$ are
listed using the notation $1^{m_1},\ldots,n^{m_n}$ which means that
the value $j$ appears with multiplicity $m_j$.  Also for brevity
$\chi(\cW,t)$ is expressed by listing its roots which are just the
exponents of $W$.  \Qed

$$\barr{c|c|c|c}
W       &\mbox{roots of $\chi(\cW,t)$}&\gamma(z)        &c_i\\[4pt]
\hline
A_n\rule{0pt}{15pt}&1,2,\ldots,n&\frac{z^{n+1}}{(1-z)^{n+1}}    &1^n\\[8pt]
B_n/C_n &1,3,\ldots,2n-1&\frac{z^{2n}(1+z)^2}{(1-z^2)^{n+1}}&1,2^{n-1}\\[8pt]
D_n     &1,3,\ldots,2n-3,n-1&\frac{z^{2n-2}(1+z)^4}{(1-z^2)^{n+1}}&1^3,2^{n-3}\\[8pt]
E_6     &1,4,5,7,8,11   &\frac{z^{12}}{(1-z)^3(1-z^2)^3(1-z^3)}&1^2,2^3,3\\[8pt]
E_7 &1,5,7,9,11,13,17&\frac{z^{18}}{(1-z)^2(1-z^2)^3(1-z^3)^2(1-z^4)}&1,2^3,3^2,4\\[8pt]
E_8     &1,7,11,13,17,19,23,29
&\frac{z^{30}}{(1-z)(1-z^2)^3(1-z^3)^2(1-z^4)^2(1-z^5)(1-z^6)}&2^2\!,3^2\!,4^2\!,5,6\\[8pt]
F_4     &1,5,7,11       &\frac{z^{12}}{(1-z)(1-z^2)^2(1-z^3)(1-z^4)}&2^2,3,4\\[8pt]
G_2     &1,5            &\frac{z^6}{(1-z)(1-z^2)(1-z^3)}&2,3
\earr$$

We should mention that Haiman~\cite[\S 7.4]{hai:cqr} independently discovered
this theorem and gave a proof which is more uniform but less elementary.
Very recently Christos Athanasiadis~\cite{ath:cps} has given
another uniform demonstration.
His main tool is the following result of Crapo and Rota~\cite{cr:cg}
which is similar 
in statement and proof to Theorem~\ref{thb} but
replaces $[-s,s]^n$ by $\bbF_p^n$ where $\bbF_p$ is the finite
field with $p$ elements, $p$ prime.
\bth[Crapo and Rota]  Let $\cA$ be any subspace arrangement in
$\bbR^n$ defined over the integers and hence over $\bbF_p$.  Then for
large enough primes $p$ we have
$$
\chi(\cA,p)=|\bbF_p^n\setm\bigcup\cA|.\quad\qed
$$
\eth
It is interesting to note that this result can also be obtained from
results of Lehrer~\cite{leh:lch} about the $l$-adic cohomology of hyperplane
complements in $\bbC^n$.  In fact Lehrer has an $l$-adic cohomological
interpretation of the 
characteristic polynomial in the equivariant case.  This suggests the
problem of trying to find versions of our two main theorems  when there
is an automorphism $g$ of $\bbC^n$ stabilizing $\cA$ and one considers the
poset of all elements of $L(\cA)$ fixed by $g$.

{\it Acknowledgment.}  We would like to thank John Stembridge who
suggested using the affine Weyl chamber to obtain the the set of
points counted in Theorem~\ref{thc}.  In addition, we thank Christos
Athanasiadis and Arun Ram for interesting discussions and relevant
references as well as the referees for helpful suggestions.

\begin{\bib}{99}

\bibitem{ath:cps} C. A. Athanasiadis, Characteristic polynomials of
subspace arrangements and finite fields, {\it Adv.\ in Math.\ } 
{\bf 122} (1996), 193--233.

\bibitem{bjo:sa}  A. Bj\"orner, Subspace arrangements, in ``Proc. 1st
European Congress Math. (Paris 1992),'' A. Joseph and R. Rentschler
eds., {\it Progress in Math.\ }, Vol.\ 122, Birkh\"auser, Boston, MA,
(1994), 321--370.

\bibitem{bl:ldt}  A. Bj\"orner and L. Lov\'asz,
Linear decision trees, subspace arrangements and M\"obius functions,
{\it \jams} {\bf 7} (1994), 667--706.

\bibitem{bly:ldt}  A. Bj\"orner, L. Lov\'asz and A. Yao,
Linear decision trees: volume estimates and topological bounds, in
``Proc. 24th ACM Symp. on Theory of Computing,''  ACM Press, New York,
NY, 1992, 170--177.

\bibitem{bs:sat}  A. Bj\"orner and B. Sagan, Subspace arrangements of
type $B_n$ and $D_n$, {\it \jac}, submitted.

\bibitem{bw:snc1}  A. Bj\"orner and M. Wachs, Shellable nonpure
complexes and posets I, {\it Trans.\ Amer.\ Math.\ Soc.\ } {\bf 348}
(1996), 1299--1327.

\bibitem{bw:snc2}  A. Bj\"orner and M. Wachs, Shellable nonpure
complexes and posets II, {\it Trans.\ Amer.\ Math.\ Soc.\ }, to appear.

\bibitem{bw:hkm} A. Bj\"orner and V. Welker, The homology of
``$k$-equal'' manifolds and related partition lattices, {\it \aim}
{\bf 110} (1995), 277--313.

\bibitem{cr:cg}   H. Crapo and G.-C. Rota, ``On the Foundations of
Combinatorial Theory: Combinatorial Geometries,''  M.I.T. Press,
Cambridge, MA, 1970.

\bibitem{hai:cqr} M. Haiman, Conjectures on the quotient ring of
diagonal invariants. {\it J. Alg.\ Combin.\ }, {\bf 3} (1994),
17--76.

\bibitem{jh:rgcg} J. Humphreys, ``Reflection Groups and Coxeter
Groups,'' Cambridge Univ.\ Press, Cambridge, England, 1990.

\bibitem{leh:lch} G. Lehrer, The $l$-adic cohomology of hyperplane
complements, {\it Bull. London Math. Soc.} {\bf 24} (1992), 76--82.

\bibitem{lin:prb}  S. Linusson, Partitions with restricted block
sizes, M\"obius functions and the $k$-of-each problem, 
{\it SIAM J. Discrete Math.\ }, to appear.

\bibitem{os:ca} P. Orlik and L. Solomon, Coxeter arrangements,
in Proc. Symp. Pure Math., Vol. 40, part 2, Amer. Math. Soc.,
Providence, RI, 1983, 269--291.

\bibitem{ot:ah} P. Orlik and H. Terao, ``Arrangements of Hyperplanes,''
Grundlehren 300, Springer-Verlag, New York, NY, 1992.

\bibitem{rot:tmf} G.-C. Rota, On the foundations of combinatorial
theory I. Theory of M\"obius functions, {\it Z. 
Wahrscheinlichkeitstheorie} {\bf 2} (1964), 340--368.

\bibitem{sta:ec1} R. P. Stanley, ``Enumerative Combinatorics,
Volume 1,''  Wadsworth and Brooks/Cole, Pacific Grove, CA, 1986. 

\bibitem{sta:gcr} R. P. Stanley, {\bf G}raph Colorings {\bf a}nd
{\bf r}elated {\bf s}ymmetric functions: {\bf i}deas and 
{\bf a}pplications, {\it Discrete Math.\ }, to appear.

\bibitem{sun:aht} S. Sundaram, Applications of the Hopf trace formula
to computing homology representations, {\it Contemp.\ Math.} {\bf 178}
(1994), 277--309.

\bibitem{sw:hrk} S. Sundaram and M. Wachs, The homology
representations of the $k$-equal partition lattice,
{\it Trans.\ Amer.\ Math.\ Soc.\ }, to appear.

\bibitem{sw:gal} S. Sundaram and V. Welker, Group actions on linear
subspace arrangements and applications to configuration spaces,
{\it Trans.\ Amer.\ Math.\ Soc.\ }, to appear.

\bibitem{ter:gef} H. Terao, Generalized exponents of a free
arrangement of hyperplanes and the Shepherd-Todd-Brieskorn formula,
{\it \im} {\bf 63} (1981), 159--179.

\bibitem{zas:grs} T. Zaslavsky, The geometry of root systems
and signed graphs, {\it Amer. Math. Monthly} {\bf 88} (1981), 88--105.

\bibitem{zas:sgc} T. Zaslavsky, Signed graph coloring, {\it \dm}  
{\bf 39} (1982) 215--228.

\bibitem{zas:cis} T. Zaslavsky, Chromatic invariants of signed graphs,
{\it \dm}  {\bf 42} (1982) 287--312.

\bibitem{zha:sba}  P. Zhang, ``Subposets of Boolean Algebras,''  Ph.D.
thesis, Michigan State University, 1994.

\end{\bib}

\end{document}